\documentclass{elsarticle}
\usepackage{amssymb}
\usepackage{amsmath}
\usepackage{amsthm}
\usepackage{setspace}
\usepackage{enumitem}
\usepackage{mathtools}
\usepackage{mathrsfs}
\usepackage{pgfplots}
\usepackage[margin=1.5 in]{geometry}

\usepackage{pgf,tikz}
\usepackage{tikz-cd}
\usepackage{mathrsfs}
\usetikzlibrary{arrows}

\newtheorem{theorem}{Theorem}[section]
\newtheorem{lemma}[theorem]{Lemma}
\newtheorem{corollary}[theorem]{Corollary}

\theoremstyle{definition}

\newtheorem{example}[theorem]{Example}

\newtheorem{conjecture}[theorem]{Conjecture}

\theoremstyle{remark}

\numberwithin{equation}{section}

\begin{document}

\begin{abstract}

	Let $N_n(F)$ denote the ring of strictly upper triangular matrices with entries in a field $F$ of characteristic zero and center $Z(N_n(F))$. We characterize the $2$-power commuting maps over $N_n(F)$, maps satisfying the identity $[f(X),X^2]=0$ for all $X\in N_n(F)$. As a consequence, we also obtain a characterization of the maps centralizing maps over $N_n(F)$, maps satisfying $[f(X),X]\in Z(N_n(F))$ for all $X\in N_n(F)$. 
	
\end{abstract}


\begin{frontmatter}

\title{Power commuting and centralizing maps on the ring of strictly upper triangular matrices}

\author{Jordan Bounds}
\address{Jordan Bounds, Department of Mathematics, Furman University, Greenville, South Carolina;}
\ead{jordan.bounds@furman.edu}

\end{frontmatter}

\section{Introduction}

Let $R$ be a ring with center $Z(R)$. A map $f:R\rightarrow R$ is called {\em centralizing} if $[f(x),x]\in Z(R)$ for all $x\in R$. If $[f(x),x]=0$ for all $x\in R$, then we say $f$ is {\em commuting}. Here we use $[a,b]=ab-ba$ to denote the standard ring commutator. Investigations into commuting maps were initiated by Posner \cite{posner57} in 1957 when he proved that a noncommutative prime ring cannot exhibit a nonzero commuting derivation. Bre{\u s}ar established the first general result regarding commuting maps, proving that an additive commuting map $f:R\rightarrow R$ where $R$ is a simple unital ring is of the form \[f(x)=\lambda x+\mu(x)\] for some $\lambda\in Z(R)$ and additive $\mu:R\rightarrow Z(R)$. Similar results regarding commuting maps satisfying particular conditions have also been established in the settings of various rings and algebras (see, for example, \cite{bflw,bresar91,bresar93,bresar95,lanski88,lanski93,vukman90}), with many concluding that commuting maps tend to be of the standard form established by Bre{\u s}ar. An overview of these results along with details on the general theory of commuting maps can be found in the survey paper \cite{bresar04} by Bre{\u s}ar.

The structures of commuting maps over various matrix rings and algebras have been explored, yielding some notable results in the cases of upper and strictly upper triangular matrices (\cite{bbc,bounds16, cheung01, franca12, franca13, franca19}). One such result, due to Beidar, Bre{\u s}ar, and Chebotar \cite{bbc}, shows that a linear commuting map defined on $T_n(F)$, the algebra of $n\times n$ upper triangular matrices with entries in a field $F$, must be of the standard form. More recently, the author \cite{bounds16} proved that linear commuting maps defined on $N_n(F)$, the ring of strictly upper triangular matrices with entries in a field $F$ of characteristic zero, are almost of the standard form. 

\begin{theorem}[Theorem 1 in \cite{bounds16}]\label{old}
Let $F$ be a field of characteristic zero, $n\ge 4$ an integer, and $N_n(F)$ the ring of strictly upper triangular matrices with entries in $F$. If $f:N_n(F)\rightarrow N_n(F)$ is a linear map satisfying $[f(X),X]=0$ for all $X\in N_n(F)$, then there exist $\lambda\in F$ and additive $\mu:N_n(F)\rightarrow\Omega$ such that \[f(X)=\lambda X+\mu(X)\] for all $X\in N_n(F)$ where $\Omega=\{ae_{1,n-1}+be_{1,n}+ce_{2,n}:a,b,c\in F\}$.
\end{theorem} 
\noindent This theorem has been recently generalized by Ko and Liu \cite{koliu23} to $N_n(R)$ where $R$ is a unital ring. Similar results have been established for the ring of infinite strictly upper triangular matrices as well \cite{slowik}.

In the present paper we seek to explore a class of related maps known as power commuting maps. We say that a map $f:R\rightarrow R$ is {\em $m$-power commuting} for a positive integer $m$ if $[f(x),x^m]=0$ for all $x\in R$. Note every commuting map is 1-power commuting. It is known that all $m$-commuting maps over certain classes of rings are also commuting (for example \cite{bflw, bresar98}). Recently, Ahmed conducted an exploration of the $m$-power commuting maps over $T_n(R)$ where $R$ is a commutative ring whose characteristic is not a divisor of $m$, concluding that such maps must be commuting \cite{ahmed19}. Ahmed and S\l owik then explored $m$-power commuting maps in the settings of $T_{\infty}(F)$ and $N_{\infty}(F)$, the rings of infinite dimensional upper and strictly upper triangular matrices, obtaining a similar classification \cite{slowik}.

In this short note we explore the $2$-power commuting maps over $N_n(F)$ where $F$ is a field of characteristic zero. Our primary result is the following.

\begin{theorem}\label{power}
Let $n\ge 4$ be an integer and $F$ a field of characteristic zero. Suppose $f:N_n(F)\rightarrow N_n(F)$ is a linear map satisfying $[f(X),X^2]=0$ for all $X\in N_n(F)$. Then there exists $\lambda\in F$ and additive $\mu:N_n\rightarrow \Psi$ such that \[f(X)=\lambda X+\mu(X)\] for every $X\in N_n(F)$ where $\Psi=\{a_1e_{1,n-1}+a_2e_{1,n}+a_3e_{2,n-1}+a_4e_{2,n}:a_i\in F\}$. 
\end{theorem}

Establishing Theorem \ref{power} has some interesting consequences. First, it shows that $2$-power commuting maps need not be commuting in this setting. Second, it allows us to obtain a complete characterization of the centralizing maps over $N_n(F)$. We 
show this explicitly by proving the following corollary.

\begin{corollary}\label{central}
Let $F$ be a field of characteristic zero and $n\ge 4$ an integer. If $g:N_n(F)\rightarrow N_n(F)$ is a linear map such that $[g(X),X]\in Z(N_n(F))$ for all $X\in N_n(F)$, then there exist $\lambda\in F$ and additive $\mu:N_n\rightarrow\Omega$ such that \[g(X)=\lambda X+\mu(X)\] for all $X\in N_n(F)$.
\end{corollary}

Throughout the remainder of this paper we will assume $n\ge 4$ is an integer and $F$ is a field of characteristic zero. We use $N_n=N_n(F)$ to denote the ring of $n\times n$ strictly upper triangular matrices with entries in $F$ and its center as $Z(N_n)$. It is known that $Z(N_n)=\{ae_{1,n}:a\in F\}$ and the identity $A^n=0$ holds for all $A\in N_n$. We will use $e_{i,j}$ to denote the standard matrix unit of $N_n$ with a 1 in position $(i,j)$ and all other entries equal to 0. For convenience, we set $\mathcal I=\{(1,n-1),(1,n),(2,n-1),(2,n)\}$ and define the sets $\Omega$ and $\Psi$ to be

\begin{align*}
\Psi&=\{A=(a_{i,j})\in N_n:a_{i,j}=0\ \text{if}\ (i,j)\notin\mathcal I\}\\
\Omega&=\{A=(a_{i,j})\in \Psi:a_{2,n-1}=0\}\\
\end{align*}

Our discussion begins with an exploration of the $2$-power commuting maps on $N_n$ in Section \ref{main}. We then explore the centralizing maps over $N_n$ in Section \ref{center}.

\section{$2$-Power Commuting Maps on $N_n$}\label{main}

We begin by providing some examples of various $2$-power commuting maps. 

\begin{example}
As has been previously stated, every commuting map is $1$-power commuting. Since \[[f(x),x^2]=x[f(x),x]+[f(x),x]x,\] it follows that every $1$ power commuting map is also $2$-power commuting. Thus, the map given in Example 2 of \cite{bounds16} should be $2$-power commuting. Indeed, if we define $g:N_n\rightarrow N_n$ by \[g(X)=\begin{pmatrix}0&\cdots&0&x_{1,2}&x_{1,n}\\ &0&\cdots&0&x_{n-1,n}\\ &&0&&0\\ &&&\ddots&\vdots\\&&&&0\end{pmatrix}\] we get $g(X)X^2=X^2g(X)=0$ for all $X\in N_n$, thus $[g(X),X^2]=0$.
\end{example}

\begin{example}
While the previous example illustrates that the image of a $2$-power commuting map must somehow account for the set $\Omega$, this set alone is insufficient for describing all $2$-power commuting maps. Consider the map $h:N_n\rightarrow N_n$ defined as \[h(X)=\begin{pmatrix}0&\cdots&0&0&x_{1,2}&x_{1,n}\\ &0&\cdots&0&x_{2,n-1}&x_{n-1,n}\\ &&0&\cdots&0&0\\ &&&\ddots&\vdots&\vdots\\&&&&&\\&&&&&0\end{pmatrix}.\]
Then we once again have $h(X)X^2=X^2h(X)=0$ for all $X\in N_n$, hence $[h(X),X^2]=0$. Moreover, we note that this example shows $2$-power commuting maps need not be commuting as \[[h(X),X]=-x_{1,2}x_{2,n-1}e_{1,n-1}-(x_{1,2}x_{2,n}+x_{1,2}x_{n-1,n})e_{1,n}+x_{2,n-1}x_{n-1,n}e_{2,n}\] is not identically zero.
\end{example}

\begin{example}\label{ex}
We can easily generalize the previous example. Let $p_{i,j}:N_n\rightarrow F$, $(i,j)\in\mathcal I$, be linear maps. We define $p:N_n\rightarrow N_n$ by \[p(X)=\begin{pmatrix}0&\cdots&0&0&p_{1,2}(X)&p_{1,n}(X)\\ &0&\cdots&0&p_{2,n-1}(X)&p_{n-1,n}(X)\\ &&0&\cdots&0&0\\ &&&\ddots&\vdots&\vdots\\&&&&&\\&&&&&0\end{pmatrix}.\]
Then $[p(X),X^2]=0$ as before. 
\end{example}

We claim that every linear $2$-power commuting map $f$ on $N_n$ can be expressed as $f(X)=\lambda X+\mu(X)$ for some $\lambda\in F$ and a map $\mu$ of the form given in Example \ref{ex}. To prove this, we first establish two lemmas that will be useful in our investigation and may be of independent interest. We note that the proofs of these lemmas rely on the process of linearization. We refer those who are unfamiliar with this process to Section 1.1 of \cite{idbook} for a detailed description.

\begin{lemma}\label{l0}
Fix $1\le i\le n-3$. Let $f$ be an additive function such that $x_{i,i+2}f(X)=0$ for all $X\in N_n$. Then $f(X)=0$ for all $X\in N_n$.
\end{lemma}
\begin{proof}
As $x_{i,i+2}f(X)=0$ for all $X\in N_n$, we have $f(e_{i,i+2})=0$. A linearization of the equation $x_{i,i+2}f(X)=0$ yields \begin{equation}\label{l0e1}
x_{i,i+2}f(Y)+y_{i,i+2}f(X)=0
\end{equation}
for all $X,Y\in N_n$. Setting $Y=e_{i,i+2}$ in Equation \ref{l0e1} we obtain $f(X)=0$ for all $X\in N_n$.
\end{proof}

\begin{lemma}\label{l1}
Suppose $g:N_n\times N_n\rightarrow N_n$ is a biadditive map satisfying $x_{2,3}g(X,X)=0$ for all $X\in N_n$. Then $g(X,X)=0$ for all $X\in N_n$.
\end{lemma}
\begin{proof}
As $x_{2,3}g(X,X)=0$ for all $X\in N_n$, we have $g(e_{2,3},e_{2,3})=0$. A linearization of the equation $x_{2,3}g(X,X)=0$ is \begin{equation}\label{l1e1}
x_{2,3}(g(Y,Z)+g(Z,Y))+y_{2,3}(g(X,Z)+g(Z,X))+z_{2,3}(g(X,Y)+g(Y,X))=0\end{equation}
for all $X,Y,Z\in N_n$.
Setting $Y=Z$ in Equation \eqref{l1e1} and using the assumption that the characteristic of $F$ is zero we obtain \begin{equation}\label{l1e2}x_{2,3}g(Y,Y)+y_{2,3}g(X,Y)+y_{2,3}g(Y,X)=0.\end{equation}
Setting $Y=e_{2,3}$, it follows that \begin{equation}\label{l1e3}g(X,e_{2,3})+g(e_{2,3},X)=0\end{equation} for all $X\in N_n$. Replacing $X$ with $X+e_{2,3}$ in the identity $x_{2,3}g(X,X)=0$ yields 
\begin{align*}0&=(x_{2,3}+1)g(X+e_{2,3},X+e_{2,3})\\
&=x_{2,3}(g(X,X)+g(X,e_{2,3})+g(e_{2,3},X)+g(e_{2,3},e_{2,3}))+g(X,X)\\
&\hspace{2.6cm}+g(X,e_{2,3})+g(e_{2,3},X)+g(e_{2,3},e_{2,3})\\
&\\
&=g(X,X).
\end{align*}
\end{proof}

Suppose $f:N_n\rightarrow N_n$ is a linear map satisfying \begin{equation}\label{e1}[f(X),X^2]=0\end{equation} for all $X\in N_n$. By the linearity of $f$ we can find linear maps $f_{i,j}:N_n\rightarrow F$ such that $f(X)=(f_{i,j}(X))$. We will frequently find it convenient to examine the entries in the matrix $[f(X),X^2]$ and equate them to 0. In general, the $(i,j)$ entry of $[f(X),X^2]$ is given by the formula \begin{equation}\label{sum}\sum\limits_{i<a<b<j}f_{i,a}(X)x_{a,b}x_{b,j}-x_{i,a}x_{a,b}f_{b,j}(X).\end{equation}

We proceed by showing there exists $\lambda\in F$ such that $f_{i,j}(X)=\lambda x_{i,j}$ for all $(i,j)\notin\mathcal I$,  and $X\in N_n$. We begin with $f_{i,i+1}(X)$.

\begin{lemma}\label{l2}
There exist $\lambda_1,\lambda_2\in F$ such that $f_{2k-1,2k}(X)=\lambda_1x_{2k-1,2k}$ and $f_{2k,2k+1}(X)=\lambda_2x_{2k,2k+1}$ for all $X\in N_n$, $1\le k$.
\end{lemma}
\begin{proof}
We use induction to show there exists $\lambda_1\in F$ such that $f_{2k-1,2k}(X)=\lambda_1x_{2k-1,2k}$. Taking the $(1,4)$ entry of $[f(X),X^2]$ and setting it equal to 0 yields \begin{equation}\label{l2e1}f_{1,2}(X)x_{2,3}x_{3,4}-x_{1,2}x_{2,3}f_{3,4}(X)=0\end{equation} for all $X\in N_n$. Defining $g:N_n\times N_n\rightarrow N_n$ by $g(X,Y)=f_{1,2}(X)y_{3,4}-y_{1,2}f_{3,4}(X)$, Equation \eqref{l2e1} becomes \begin{equation}\label{l2e2}x_{2,3}g(X,X)=0\end{equation} for all $X\in N_n$. Applying Lemma \ref{l1} we conclude \begin{equation}\label{l2e3}f_{1,2}(X)x_{3,4}-x_{1,2}f_{3,4}(X)=g(X,X)=0\end{equation} for all $X\in N_n$.

Setting $X=e_{3,4}$ in Equation \eqref{l2e3} we observe $f_{1,2}(e_{3,4})=0$. Thus, replacing $X$ with $X+e_{3,4}$ in Equation \eqref{l2e3} we obtain 
\begin{align*}0&=f_{1,2}(X+e_{3,4})(x_{3,4}+1)-x_{1,2}f_{3,4}(X+e_{3,4})\\
&=f_{1,2}(X)x_{3,4}+f_{1,2}(X)-x_{1,2}f_{3,4}(X)-x_{1,2}f_{3,4}(e_{3,4})\\
&=f_{1,2}(X)-x_{1,2}f_{3,4}(e_{3,4}).
\end{align*}
Therefore, $f_{1,2}(X)=\lambda_1x_{1,2}$ where $\lambda_1=f_{3,4}(e_{3,4})$.

Suppose $f_{2k-1,2k}(X)=\lambda_1x_{2k-1,2k}$ for some $k\ge 1$. Taking the $(2k-1,2k+2)$ entry of $[f(X),X^2]$ and setting it equal to 0 yields \begin{align*}
0&=f_{2k-1,2k}(X)x_{2k+1,2k+2}-x_{2k-1,2k}f_{2k+1,2k+2}(X)\\
&=\lambda_1x_{2k-1,2k}x_{2k+1,2k+2}-x_{2k-1,2k}f_{2k+1,2k+2}(X)\\
&=x_{2k-1,2k}(\lambda_1x_{2k+1,2k+2}-f_{2k+1,2k+2}(X)).
\end{align*}

It follows from Lemma \ref{l0} that $f_{2k+1,2k+2}(X)=\lambda_1x_{2k+1,2k+2}$. Applying a similar argument using the $(2,5)$ entry of $[f(X),X^2]$ will yield $f_{2k,2k+1}(X)=\lambda_2x_{2k,2k+1}$ for some $\lambda_2\in F$.
\end{proof}

\begin{lemma}\label{l3}
There exists $\lambda\in F$ such that $f_{i,i+1}(X)=\lambda x_{i,i+1}$ for all $1\le i\le n-1$.
\end{lemma}
\begin{proof}
It suffices to show $\lambda_1$ and $\lambda_2$ from Lemma \ref{l2} are equal. Taking the $(1,5)$ entry of $[f(X),X^2]$ and setting it equal to 0 yields \begin{align*}0&=f_{1,2}(X)x_{2,3}x_{3,5}+f_{1,2}(X)x_{2,4}x_{4,5}+f_{1,3}(X)x_{3,4}x_{4,5}\\
&\hspace{2.5cm}-x_{1,2}x_{2,3}f_{3,5}(X)-x_{1,2}x_{2,4}f_{4,5}(X)-x_{1,3}x_{3,4}f_{4,5}(X)\\
&\\
&=x_{1,2}x_{2,3}(\lambda_1x_{3,5}-f_{3,5}(X))+x_{1,2}x_{2,4}x_{4,5}(\lambda_1-\lambda_2)+x_{3,4}x_{4,5}(f_{1,3}(X)-\lambda_2x_{1,3})
\end{align*} for all $X\in N_n$ where the second equality is due to Lemma \ref{l2}. Setting $X=e_{1,2}+e_{2,4}+e_{4,5}$ yields $\lambda_1=\lambda_2$.
\end{proof}

\begin{lemma}\label{l4}
There exists $\lambda\in F$ such that $f_{i,j}(X)=\lambda x_{i,j}$ for all $X\in N_n$ and $(i,j)\notin\mathcal I$.
\end{lemma}
\begin{proof}
We proceed by double induction on the integers $i,k \ge 1$ where $k=j-i$. In this manner, we show $f_{i,i+k}(X)=\lambda x_{i,i+k}$, $(i,i+k)\notin\mathcal I$, by working along the diagonals of the matrix $f(X)$.

For $i=k=1$, we have $f_{1,2}(X)=\lambda x_{1,2}$ by Lemma \ref{l3} for some $\lambda\in F$. If we suppose $f_{i,i+1}(X)=\lambda x_{i,i+1}$, then $f_{i+1,i+2}(X)=\lambda x_{i+1,i+2}$, also by Lemma \ref{l3}. So suppose for each $1\le t\le k$ we have $f_{i,i+t}(X)=\lambda x_{i,i+t}$ for all $X\in N_n$, $i\ge 1$. We show $f_{i+1,i+k+1}(X)=\lambda x_{i+1,i+k+1}$. As $k=j-i$, we will work with $f_{i+1,j+1}$ for convenience.

Taking the $(i+1,j+3)$ entry of $[f(X),X^2]$ and setting it equal to 0 we obtain \begin{equation}\label{l4e1}\sum\limits_{i+1<a<b<j+3}[f_{i+1,a}(X)x_{a,b}x_{b,j+3}-x_{i,a}x_{a,b}f_{b,j+3}(X)]=0.
\end{equation}
Applying the inductive assumption we observe
\begin{align*}
&\sum\limits_{\substack{i+1<a<b<j+3\\ (a,b)\ne (j+1,j+2)\\ (a,b)\ne(i+2,i+3)}}[f_{i+1,a}(X)x_{a,b}x_{b,j+3}-x_{i+1,a}x_{a,b}f_{b,j+3}(X)]\\
&\\
=&\sum\limits_{\substack{i+1<a<b<j+3\\ (a,b)\ne (j+1,j+2)\\ (a,b)\ne(i+2,i+3)}}[\lambda x_{i+1,a}x_{a,b}x_{b,j+3}-x_{i+1,a}x_{a,b}\lambda x_{b,j+3}]\\
&\\
=&0.
\end{align*} 
Therefore,
\begin{align*}
0&=\sum\limits_{i+1<a<b<j+3}[f_{i+1,a}(X)x_{a,b}x_{b,j+3}-x_{i+1,a}x_{a,b}f_{b,j+3}(X)]\\
&=f_{i+1,j+1}(X)x_{j+1,j+2}x_{j+2,j+3}-x_{i+1,j+1}x_{j+1,j+2}\lambda x_{j+2,j+3}\\
&\hspace{1cm} +\lambda x_{i+1,i+2}x_{i+2,i+3}x_{i+3,j+3}-x_{i+1,i+2}x_{i+2,i+3}f_{i+3,j+3}(X)\\
&\\
&=x_{j+1,j+2}x_{j+2,j+3}(f_{i+1,j+1}(X)-\lambda x_{i+1,j+1})+x_{i+1,i+2}x_{i+2,i+3}(\lambda x_{i+3,j+3}-f_{i+3,j+3}(X))\tag{\theequation}\label{l4e2}
\end{align*}
for all $X\in N_n$.

Setting $E=e_{j+1,j+2}+e_{j+2,j+3}$ we observe $f_{i+1,j+1}(E)=0$. Thus, replacing $X$ with $X+E$ in Equation \eqref{l4e2} yields

\begin{multline}\label{l4e3}
0=x_{j+1,j+2}(f_{i+1,j+1}(X)-\lambda x_{i+1,j+1})+x_{j+2,j+3}(f_{i+1,j+1}(X)-\lambda x_{i+1,j+1})\\ +f_{i+1,j+1}(X)-\lambda x_{i+1,j+1}-x_{i+1,i+2}x_{i+2,i+3}f_{i+3,j+3}(E).
\end{multline}
Substituting $X=e_{j+1,j+2}$ into Equation \eqref{l4e3} yields $f_{i+1,j+1}(e_{j+1,j+2})=0$ as $F$ has characteristic different from 2. Thus, by replacing $X$ with $X+e_{j+1,j+2}$ in Equation \eqref{l4e3} we obtain
\begin{equation}
0=f_{i+1,j+1}(X)-\lambda x_{i+1,j+1}.
\end{equation}
Therefore, $f_{i+1,j+1}(X)=\lambda x_{i+1,j+1}$ as desired.
\end{proof}

We now prove Theorem \ref{power}

\begin{proof}[Proof of Theorem \ref{power}]
Applying the linearity of $f$ we obtain linear functions $f_{i,j}:N_n\rightarrow F$ such that $f(X)=(f_{i,j}(X))$ for all $X\in N_n$. By Lemma \ref{l4} there exists $\lambda\in F$ such that $f_{i,j}(X)=\lambda x_{i,j}$ for all $(i,j)\notin\mathcal I$, $X\in N_n$. For each $(i,j)$ we define $\mu_{i,j}:N_n\rightarrow F$ to be \[\mu_{i,j}(X)=\begin{cases}f_{i,j}(X)-\lambda x_{i,j}\ \text{if}\ (i,j)\in\mathcal I\\ 0\ \hspace{2.1cm}\text{if}\ (i,j)\notin\mathcal I\end{cases}.\] Taking $\mu:N_n\rightarrow\Psi$ to be the map defined by $\mu(X)=(\mu_{i,j}(X))$ we have \[f(X)=\lambda X+\mu(X)\] for all $X\in N_n$. Given $A,B\in N_n$ we have \[\mu(A+B)=f(A+B)-\lambda(A+B)=f(A)-\lambda A+f(B)-\lambda B=\mu(A)+\mu(B),\] hence $\mu$ is additive.
\end{proof}

We believe a modification of this argument can be used to characterize the $m$-power commuting maps over $N_n$. For each $m\ge 2$ define $\mathcal I_m$ to be \[\mathcal I_m=\{(i,j):1\le i\le m, n-m+1\le j\le n\}.\] Let $\Psi_m$ be the following subset of $N_n$: \[\Psi_m=\{A=(a_{i,j})\in N_n:a_{i,j}=0\ \text{if}\ (i,j)\notin\mathcal I_m\}.\]
We believe the following statement to be true.
\begin{conjecture}
Let $1<m<n-1$ be an integer and suppose $f:N_n\rightarrow N_n$ is a linear $m$-power commuting map. Then there exist $\lambda\in F$ and additive $\mu:N_n\rightarrow\Psi_m$ such that \[f(X)=\lambda X+\mu(X).\]
\end{conjecture}
The case $m=1$ is unique in that the image of the corresponding $\mu$ map lies in the set $\Omega$, which does not follow the same structure as the sets $\Psi_m$ for $m>1$. Additionally, since $A^n=0$ for every $A\in N_n(F)$, every map is trivially $m$-power commuting for $m\ge n-1$. The above conjecture appears to be a reasonable means for describing the remaining values of $m$. 

\section{Centralizing Maps on $N_n$}\label{center}
Recall a map $f:R\rightarrow R$ is called {\em centralizing} if $[f(x),x]\in Z(R)$ for all $x\in R$. It is clear from their definitions that every commuting map is necessarily centralizing. It is also known that, under mild conditions, the converse is true in the setting of prime rings (see Proposition 3.1 in \cite{bresar93}). However, it is a simple task to construct examples of maps on $N_n$ that are centralizing but not commuting.

\begin{example}\label{not com}
Consider the map $p$ defined by \[p(A)=\begin{pmatrix}0&\cdots&0&a_{n-1,n}&a_{1,n}\\ &0&\cdots&0 &a_{1,2}\\ &&0&\cdots&0\\ &&&\ddots&\vdots\\ &&&&0\end{pmatrix}\] for all $A=(a_{i,j})\in N_n$. Then $p$ satisfies \[ [p(A),A]=(a_{n-1,n}^2-a_{1,2}^2)e_{1,n}\in Z(N_n)\] and is therefore centralizing. As the value $a_{n-1,n}^2-a_{1,2}^2$ is not identically 0, however, $p$ is not commuting.
\end{example}

Using Theorem \ref{power} we can obtain an explicit description of the centralizing maps on $N_n$ by proving Corollary \ref{central}.

\begin{proof}[Proof of Corollary \ref{central}]
We have $XZ(N_n)=Z(N_n)X=\{0\}$ for all $X\in N_n$. Thus, \[[f(X),X^2]=X[f(X),X]+[f(X),X]X=0\] for all $X\in N_n$, making $f$ a 2-power commuting map on $N_n$. Applying Theorem \ref{power} we obtain $\lambda\in F$ and additive $\mu:N_n\rightarrow\Psi$ such that \[f(X)=\lambda X+\mu(X)\] for every $X\in N_n$. It follows that \begin{align*}
[f(X),X]&=[\lambda X+\mu(X),X]\\
&=[\mu(X),X]\in Z(N_n).
\end{align*}
But $[\mu(X),X]=-\mu_{2,n-1}(X)x_{1,2}e_{1,n-1}+\mu_{2,n-1}(X)x_{n-1,n}e_{2,n}+(\mu_{1,n-1}(X)x_{n-1,n}-\mu_{2,n}(X)x_{1,2})e_{1,n}$. Since $Z(N_n)=\{ae_{1,n}:a\in F\}$, it must be the case that $\mu_{2,n-1}(X)=0$ for all $X$. Therefore, $\mu(N_n)\subset\Omega$ as desired.
\end{proof}

\section*{Acknowledgement}
This material is based upon work supported by the National Science Foundation under Award No. 2316995.

\bibliographystyle{amsplain}

\end{document}